\documentclass{elsart3-1}

\usepackage{amssymb,amsmath}

\usepackage[english,francais]{babel}
\usepackage{graphicx,color}

\newtheorem{theorem}{Theorem}[section]

\newtheorem{e-proposition}[theorem]{Proposition}

\newtheorem{e-definition}[theorem]{Definition\rm}

\renewcommand{\theequation}{\arabic{equation}}
\def\R{\mathbb{R}}

\newcommand{\beq}{\begin{equation}}
\newcommand{\eeq}{\end{equation}}
\newcommand{\fer}{\eqref}
\newcommand{\da}{\delta}
\newcommand{\f}{\frac}

\def\og{\leavevmode\raise.3ex\hbox{$\scriptscriptstyle\langle\!\langle$~}}
\def\fg{\leavevmode\raise.3ex\hbox{~$\!\scriptscriptstyle\,\rangle\!\rangle$}}
\def\p{\partial}
\def\e{\varepsilon}

\journal{the Acad\'emie des sciences}
\begin{document}
% place in the next line the header (rubrique) chosen for your article,
% if you know it (you can also have 2, format : Header1/Header2
\centerline{}
\begin{frontmatter}


% \author{Name\thanksref{label2}}
\selectlanguage{english}
\title{A Hamilton-Jacobi method to describe the evolutionary equilibria \\in heterogeneous environments and with non-vanishing effects of mutations}

% The [label1] can be suppressed if there is only one address for all authors

\selectlanguage{english}
\author[authorlabel1]{Sylvain Gandon}
\ead{sylvain.gandon@cefe.cnrs.fr}
\author[authorlabel2]{Sepideh Mirrahimi},
\ead{sepideh.mirrahimi@math.univ-toulouse.fr}

\address[authorlabel1]{ Centre d'Ecologie Fonctionnelle et Evolutive (CEFE), UMR CNRS 5175, Montpellier, France}
\address[authorlabel2]{CNRS, Institut de Math\'ematiques (UMR CNRS 5219), Universit\'e Paul Sabatier, 118 Route de Narbonne, 31062 Toulouse Cedex, France}

% If you know the dates of reception, and acceptation you can put them now;
%  idem the name of the person presenting the Note

\medskip
\begin{center}
{\small Received *****; accepted after revision +++++\\
Presented by £££££}
\end{center}

\begin{abstract}
\selectlanguage{english}
In this note, we characterize the solution of a system of  elliptic integro-differential equations describing a phenotypically structured population subject to mutation, selection and migration. Generalizing an approach based on Hamilton-Jacobi equations, we identify the dominant terms of the  solution when the mutation term is small (but nonzero). This method was initially used, for different problems from evolutionary biology, to identify the asymptotic solutions, while the mutations vanish, as a sum of Dirac masses. A key point is a uniqueness property related to the weak KAM theory. This method allows to go further than the Gaussian approximation commonly used by biologists and is an attempt to fill the gap between the theories of adaptive dynamics and
quantitative genetics.% and to make a link between the theories of adaptive dynamics and quantitative genetics.

\vskip 0.5\baselineskip

\selectlanguage{francais}
{\bf R\'esum\'e} \vskip 0.5\baselineskip 
{\bf M\'ethode Hamilton-Jacobi pour d\'ecrire des \'equilibres \'evolutives dans les environnements h\'et\'erog\`enes et avec des mutations non-\'evanescentes.}
Dans cette note, nous \'etudions un syst\`eme d'\'equations int\'egro-diff\'erentielles elliptiques, d\'ecrivant une population structur\'ee par trait ph\'enotypique soumise \`a des mutations, \`a la s\'election et \`a des migrations. Nous g\'en\'eralisons une approche bas\'ee sur des \'equations de Hamilton-Jacobi pour d\'et\'erminer les termes dominants de la solution lorsque les effets des mutations sont petits (mais non-nuls). Cette m\'ethode \'etait initialement utilis\'ee, pour diff\'erents probl\`emes venant de la biologie \'evolutive, pour identifier les solutions asymptotiques, lorsque les effets des mutations tendent vers $0$, sous forme des sommes de masses de Dirac. Un point cl\'e est une propri\'et\'e d'unicit\'e en rapport avec la th\'eorie de KAM faible. Cette m\'ethode nous permet d'aller au-del\`a des approximations Gaussiennes habituellement utilis\'ees par les biologistes et contribue ainsi \`a relier les th\'eories de la dynamique adaptative et de la
g\'en\'etique quantitative. 

\end{abstract}
\end{frontmatter}

%
%% now the Version française abrégée, if it exists
%\selectlanguage{francais}
%\section*{Version fran\c{c}aise abr\'eg\'ee}
%% Text of your Version française abrégée here.
%% Note you do not need to repeat here equations that you use in the
%% main text - for example 'voir (3)' is quite acceptable.
%
%

\selectlanguage{english}
% main text
\section{Introduction}

During the last decade, an approach based on Hamilton-Jacobi equations with constraints has been developed to describe the asymptotic evolutionary dynamics of phenotypically structured populations, in the limit of vanishing mutations. The mathematical modeling of such phenomena leads to parabolic (or elliptic for the steady case) integro-differential equations, whose solutions tend, as the diffusion term vanishes, to a sum of Dirac masses, corresponding to dominant traits. These asymptotic solutions can be described using the Hamilton-Jacobi approach. There is a large literature on this method. We refer to \cite{OD.PJ.SM.BP:05,GB.BP:08,SM:11} for the establishment of the basis of this approach for problems from evolutionary biology. Note that related tools were already used in the case of local equations (for instance KPP type equations) to describe the propagation phenomena (see for instance \cite{MF:86,LE.PS:89}).

In almost all the previous works, the Hamilton-Jacobi approach has been used to describe the limit of the solution, corresponding to the population's phenotypical distribution,  as the mutations' steps vanish. However, from the biological point of view, it is sometimes more relevant to consider non-vanishing mutation's steps. A recent work \cite{SM.JR:15-1} has pointed out that such tools can also be used, for a simple model with  homogeneous environment, to characterize the solution while mutations' steps are small but nonzero. In this note, we show how such  results can be obtained in a more complex situation with a heterogeneous environment.
\\

Our purpose in this note is to study the  solutions to the following system, for $z\in  \R$,
\begin{equation}
\label{main}
\begin{cases}
-\e^2  n_{\e,1}''(z)=n_{\e,1}(z) R_1(z,N_{\e,1})+m_2n_{\e,2}(z)-m_1n_{\e,1}(z),\\
-\e^2   n_{\e,2}''(z)=n_{\e,2}(z) R_2(z,N_{\e,2})+m_1n_{\e,1}(z)-m_2n_{\e,2}(z),\\
N_{\e,i}=\int_\R n_{\e,i}(z)dz, \qquad \text{for $i=1,2$},
\end{cases}
\end{equation}
with 
%\begin{equation}
%\label{Ni}
%N_{\e,i}=\int_\R n_{\e,i}(z)dz, \qquad \text{for $i=1,2$},
%\end{equation}
%and 
\beq
\label{Ri}
R_i(z,N_i)= r_i- g_i(z-\theta_i)^2-\kappa_i N_i,\qquad \text{with $\theta_1=-\theta$ and $\theta_2=\theta$}.
\eeq
This system represents the equilibrium of a phenotypically structured population under mutation, selection and migration between two habitats. For more details on the modeling and the biological motivations, see Section \ref{sec:mod}.

Note that the asymptotic behavior, as $\e\to 0$ and along subsequences, of the solutions to this system, under the assumption $m_i>0$, for $i=1,2$, and for bounded domains, was already studied in \cite{SM:12}. In the present work, we go further than the asymptotic limit along subsequences and we obtain uniqueness of the limit and identify the dominant terms of the solution when $\e$ is small but nonzero.  
\\

 {\bf The main elements of the method:}\\
To describe the solutions $n_{\e,i}(z)$ we use a WKB ansatz
$$
n_{\e,i}(z)=\frac{1}{\sqrt{2\pi \e}} \exp \left(\frac{u_{\e,i}(z)}{\e} \right).
$$
Note that a first approximation that is commonly used in the theory of 'quantitative genetics' (a theory in evolutionary biology that investigates the evolution of continuously varying traits \cite{Rice-book}-chapter $7$), is a gaussian distribution of the following form
$$
n_{\e,i}(z)=\frac{N_i}{\sqrt{2\pi\e}\sigma}\exp  \left(\frac{-(z-z^*)^2}{\e\sigma^2} \right)
=\frac{1}{\sqrt{2\pi\e}} \exp \left(\frac{-\frac{1}{2\sigma^2}(z-z^*)^2+\e\log \frac{N_i}{\sigma}}{\e} \right).
$$
Here, we try to go further than this a priori gaussian assumption and to approximate directly $u_{\e,i}$. To this end, we write an expansion for 
$u_{\e,i}$ in terms of $\e$:
\beq
\label{ap-ue}
u_{\e,i}= u_i+\e v_i+\e^2 w_i+O(\e^3).
\eeq
 We prove that   $u_1=u_2=u$ is the unique viscosity solution to a Hamilton-Jacobi equation with constraint. The uniqueness of the viscosity solution to such Hamilton-Jacobi equation with constraint is related to  the uniqueness of the Evolutionary Stable Strategy (ESS), see Section \ref{sec:ad} for definition  and for the result on the uniqueness of the ESS, and  to the weak KAM theory \cite{AF:16}.
In section \ref{sec:u},  we compute explicitly  $u$ which indeed satisfies 
$$
\max_\R \; u(z)= 0,
$$
with the maximum points attained at one or two points corresponding to the ESS points of the problem.
We then notice that, while $u(z)<0$, $n_{\e,i}(z)$ is exponentially small. Therefore, only the values of $v_i$ and $w_i$ close to the zero level sets of $u$ matter, i.e. the ESS points. In section \ref{sec:vw}, we provide the main elements to compute formally $v_i$ and hence its second order Taylor expansion around the ESS points  and the value of $w_i$ at those points. Then, we show, in section \ref{sec:mom}, that these approximations together with  a fourth order Taylor expansion of $u$ around the ESS points are enough to approximate the moments of the population's distribution with an error of order $\e^2$.

The mathematical details of our results will be provided in \cite{SM:17}. The biological applications will be detailed in \cite{SG.SM:17}.

\section{Model and motivation}
\label{sec:mod}

The solution of \eqref{main} corresponds to the steady solution to the following system, for $(t,z)\in \R^+\times \R$,
\begin{equation}
\label{dyn}
\begin{cases}
\p_t n_i(t,z)-\e^2 \f{\p^2}{\p z^2} n_i(t,z)=n_i(t,z) R_i(z,N_i(t))+m_jn_j(t,z)-m_in_i(t,z),\quad i=1,2,\quad j=2,1,\\
N_i(t)=\int_\R n_i(t,z)dz, \qquad \text{for $i=1,2$}.
\end{cases}
\end{equation}
%with 
%\begin{equation}
%\label{Ni}
%N_i(t)=\int_\R n_i(t,z)dz, \qquad \text{for $i=1,2$}.
%\end{equation} 
This system represents the dynamics of a population that is structured by a phenotypical trait $z$, and which lives in two habitats. We denote by $n_i(t,z)$ the density of the phenotypical distribution in habitat $i$, and by $N_i$ the total population's size in habitat $i$. The growth rate $R_i(z,N_i)$ is given by \fer{Ri}, where $r_i$ represents the maximum intrinsic growth rate, $g_i$ is the strength of the selection, $\theta_i$ is the optimal trait in habitat $i$ and $\kappa_i$ represents the intensity of the competition. The constants $m_i$ are the migration rates between the habitats. In this note we assume that there is positive migration rate in both directions, i.e.
\begin{equation}
\label{as:m}
m_i>0, \qquad \text{i=1,2}.
\end{equation}
However, the source and sink case, where for instance $m_2=0$, can also be analyzed using similar tools. We refer to \cite{SG.SM:17} for the analysis of this case. We additionally assume that
\begin{equation}
\label{as:r-m}
\max( r_1-m_1,r_2-m_2) > 0.
\end{equation}
This guarantees that the population does not get extinct.
\\

Such phenomena have already  been studied by several approaches. A first class of results  are based on the adaptive dynamics approach, where one considers that the mutations are very rare such that the population has time to attain its equilibrium between two mutations and hence the population's distribution has discrete support (one or two points in a two habitats model) \cite{GM.IC.SG:97,TD:00,CF.SM.EM:12}.  A second class of results are based on an approach known as  'quantitative genetics', which allows more frequent mutations and does not  separate the evolutionary and the ecological time scales . A main assumption in this class of works is that one considers that the population's distribution is a gaussian \cite{AH.TD.ET:01,OR.MK:01} or, to take into account the possibility of dimorphic populations, a sum of one or two gaussian distributions \cite{SY.FG:09,FD.OR.SG:13}. 

In our work, as in the  quantitative genetics framework, we also consider continuous phenotypical distributions. However,  we don't assume any a priori gaussian assumption. We compute directly the population's distribution and in this way we correct the previous approximations.  To this end, we also provide some results in the framework of adaptive dynamics and in particular, we generalize previous results on the identification of the ESS to the case of nonsymetric habitats. Furthermore, our work makes a connection between the two approaches of adaptive dynamics and quantitative genetics.

\section{The adaptive dynamics framework}
\label{sec:ad}

In this section, we introduce some notions  from the theory of adaptive dynamics  that we will be using in the next sections  \cite{GM.IC.SG:97}.  We also provide our main result in this framework.

{\bf Effective fitness:} 
 The effective fitness $W(z;N_1,N_2)$ is the largest eigenvalue of the following matrix:
\beq
\label{efitness}
\mathcal A (z;N_1,N_2)= \left( 
\begin{array}{cc}
R_1(z ; N_1) -m_1 & m_2\\
m_1 & R_2(z ; N_2) -m_2
\end{array}
\right)
\eeq
This indeed corresponds to the \emph{effective} growth rate associated with trait $z$ in the whole metapopulation when the total population's sizes are given by $(N_1,N_2)$.

{\bf Demographic equilibrium:} 
Consider a set of points $\Omega=\{z_1,\cdots z_m\}$. The demographic equilibrium corresponding to this set is given by $(n_1(z),n_2(z))$, with the total population's  sizes $(N_1,N_2)$, such that
$$
n_i(z)=\sum_{j=1}^m \alpha_{i,j}\delta(z-z_j),\qquad N_i = \sum_{j=1}^m \alpha_{i,j},\qquad W(z_j,N_1,N_2)=0,
$$
and such that $(\alpha_{1,j},\alpha_{2,j})^{T}$  is the right eigenvector associated with the largest eigenvalue $W(z_j,N_1,N_2)=0$ of $\mathcal A(z_j;N_1,N_2)$.

{\bf Evolutionary stable strategy:} 
A set of points $\Omega^*=\{z_1^*,\cdots, z_m^*\}$ is called an evolutionary stable strategy (ESS) if
$$
W(z,N^*_1,N_2^*)=0,\quad \text{for $z\in \mathcal A$ and} ,\quad W(z,N_1^*,N_2^*)\leq 0, \quad \text{for  $z\not\in \mathcal A$,}
$$
where $(N_1^*,N_2^*)$ is the total population's sizes corresponding to the demographic equilibrium associated with the set $\Omega^*$.

Since, there are only two habitats we expect that at most two distinct traits coexist at the evolutionary stable equilibrium. We prove indeed the following:

\begin{theorem}
Assume \fer{as:m}--\fer{as:r-m}.   There exists a unique set of points $\Omega^*$ which is an evolutionary stable strategy. Such set has at most two elements.
\end{theorem}

 We call an evolutionary stable strategy which has one (respectively two) element, a monomorphic (respectively dimorphic) ESS. We can indeed give a criterion to have monomorphic or dimorphic ESS, and we can identify the dimorphic ESS in the general case (see \cite{SM:17} for more details).

\section{How to compute the zero order terms $u_i$}
\label{sec:u}
The identification of the zero order terms $u_i$ is based on the following result. Note that the  part (ii) of the  theorem below is a variant of Theorem 1.1 in \cite{SM:12}. 

\begin{theorem}\label{thm:main}
Assume \fer{as:m}--\fer{as:r-m}. \\
(i) As $\e \to 0$, $(n_{\e,1}, n_{\e,2})$ converges to $(n_1^*,n_2^*)$, the demographic equilibrium of the unique  ESS of the model. Moreover, as $\e \to 0$, $N_{\e,i}$  converges to $N_i^*$, the total population's size in patch $i$ corresponding to this demographic equilibrium.
\\
(ii) As $\e\to 0$, both sequences $(u_{\e,i})_\e$, for $i=1,2$, converge along subsequences and locally uniformly in $\R$ to a continuous function $u\in \mathrm{C}(\R)$,  such that $u$ is a viscosity solution to the following equation
\begin{equation}
\label{HJ}
\left\{ \begin{array}{ll}-|u'|^2= W(z,N_1^*,N_2^*),&\quad \text{in $\R$},\\
\max_{z\in \R}u(z)=0.\end{array}\right.
\end{equation}
Moreover, we have the following condition on  the zero level set of $u$:
$$
{\rm supp}\, n_1^*= {\rm supp}\, n_2^*\subset \{z\, |\, u(z)=0 \} \subset \{z\, |\, W(z,N_1^*,N_2^*)=0\}.
$$
(iii) There exists constants $(\lambda_i,\nu_i)$, for $i=1,2$, which can be determined explicitly  from $m_1$, $m_2$, $g_1$, $g_2$, $\kappa_1$, $\kappa_2$ and $\theta$, such that, under the condition 
\beq
\label{non-deg}
r_2 \neq \lambda_1r_1+\nu_1,\qquad r_1\neq  \lambda_2r_2+\nu_2,
\eeq
we have
\beq
\label{aubry}
{\rm supp}\, n_1^*= {\rm supp}\, n_2^*= \{z\, |\, u(z)=0 \} = \{z\, |\, W(z,N_1^*,N_2^*)=0\}.
\eeq
 The solution of \eqref{HJ}--\eqref{aubry} is unique and hence the whole sequence $(u_{\e,i})_\e$   converge  locally uniformly in $\R$ to   $u$.
\end{theorem}

Note that a Hamilton-Jacobi equation of type \fer{HJ} in general might admit several viscosity solutions. Here, the uniqueness is obtained thanks to \fer{aubry} and a property from   the weak KAM theory, which is the fact that the viscosity solutions are completely determined by one value taken on each static class of the Aubry set (\cite{PL:82}, Chapter 5 and \cite{GC:01}).
 In what follows we assume that \fer{non-deg} and hence \fer{aubry} always hold. We then give an explicit formula for $u$ considering two cases:

(i) {\bf Monomorphic ESS : } We consider the case where there exists a unique monomorphic ESS $z^{*}$ and the corresponding demographic equilibrium is given by $(N_1^{*}\da(z^{*}),N_2^{*}\da(z^{*}))$. Then  $u$ is given by
\beq
\label{u-exp}
u(z)= -\big| \int_{z^{*}}^{ z} \sqrt{- W(x; N_1^{*}, N_2^{*})} dx \big|.
\eeq

(ii) {\bf Dimorphic ESS : } We next consider the case where there exists a unique dimorphic ESS $(z_a^{*},z_b^{*})$   with the demographic equilibrium:
$
n_i=\nu_{a,i} \da(z-z^*_a)+\nu_{b,i} \da(z-z^*_b),
$ and
$ \nu_{a,i}+ \nu_{b,i} =N_i^*$. Then  $u$ is given by
$$
 u(z)=\max \Big( - |\int_{z_a^{*}}^{ z} \sqrt{- W(x; N_1^{*}, N_2^{*})} dx|
  , - |\int_{z_b^{*}}^{z} \sqrt{- W(x; N_1^{*}, N_2^{*})} dx |\Big).
$$

\section{How to compute the next order terms}
\label{sec:vw}

%
%We notice that away from the ESS points the zero order term $u$ is strictly negative and hence $n_{\e,i}$ is exponentially small. Therefore, in the characterization of the population's distribution, only the values of $v_i$ and $w_i$ around the ESS points matter.  It is indeed enough to compute a fourth order approximation of $u$ and a second order approximation of $v_i$ around the ESS points and to determine the value of $w_i$ at those points, in order to approximate the moments of the distribution with an error of order $\e^2$. In the previous section, we computed explicitly $u$ and hence one can easily compute also its  fourth order Taylor expansion.

  In this section, we give the main elements to compute formally $v_i$ and the value of $w_i$ at the ESS point, with $v_i$ and $w_i$ the correctors introduced by \fer{ap-ue},   in the case of monomorphic population.  For the details of the computations for both monomorphic and dimorphic populations, we refer the interested reader to \cite{SM:17}.

We consider  the case of monomorphic   population where the demographic equilibrium corresponding to the monomorphic ESS  is given by $(N_1^{*}\da(z-z^{*}), N_2^{*}\da(z-z^{*}) )$. One can compute, using \fer{u-exp}, a Taylor expansion of order $4$ around the ESS point $z^{*}$:
$$
u(z) = -\f{A}{2 }(z-z^{*})^2+B(z-z^{*})^3+C(z-z^{*})^4+O(z-z^{*})^5.
$$
To provide an approximation of the moments of the population's distribution, we have to compute constants $D_i$, $E_i$  and $F_i$ such that
$$
v_i(z)=v_i(z^{*})+D_i (z-z^{*}) +E_i(z-z^{*})^2+O(z-z^{*})^3,\qquad w_i(z^{*})=F_i.
$$
A first element of the computations is obtained by replacing the functions $u$, $v_i$ and $w_i$ by the above approximations to compute $N_{\e,i}=\int_\R  n_{\e,i} (z)dz$. This leads to
$$
v_i(z^{*})= \log \big(N_i^{*} \sqrt{A} \big),
\quad
N_{\e,i} =N_i^{*} +\e K_i^{*}+O(\e^2),\quad \text{with }\quad K_i^{*}= N_i^{*} \big(  \f{3C}{A^2}+\f{E_i}{A}+F_i  \big).
$$
Note also that writing \fer{main} in terms of $u_{\e,i}$ we obtain
\begin{equation}
\label{main2}
\begin{cases}
-\e  u_{\e,1}''(z) = |u_{\e,1}'|^2+ R_1(z,N_{\e,1})+m_2\exp \big( \f{u_{\e,2}-u_{\e,1}}{\e}\big)-m_1 ,\\
-\e  u_{\e,2}''(z) = |u_{\e,2}'|^2+ R_2(z,N_{\e,2})+m_1\exp \big( \f{u_{\e,1}-u_{\e,2}}{\e}\big)-  m_2.
\end{cases}
\end{equation}
A second element, is obtained by keeping the zero order terms in  the first line of \fer{main2} and using \fer{HJ} to obtain 
\beq
\label{Q2/Q1}
  v_2(z)-v_1(z) = \log \Big( \f{1}{m_2}\big( W(z,N_1^*,N_2^*)- R_1(z,N_1^*)+m_1 \big) \Big).
\eeq
The last element, is derived from keeping the terms of order $\e$ in \fer{main2} which leads to
\beq
\label{1-order}
-  u''=2  u'   v_i'-\kappa_i K_i ^{*}+m_j \exp(v_j-v_i)(w_j-w_i),\quad \text{for $\{i,j\} =\{1,2\}$}.
\eeq
The functions $v_i$ and the coefficients $D_i$, $E_i$ and $F_i$ can be computed by combining the above elements.

\section{Approximation of the moments}
\label{sec:mom}

The above approximations of $u$, $v_i$ and $w_i$ around the ESS points allow us to estimate the moments of the population's distribution. In the monomorphic case these approximations are given below:
$$
\begin{cases}
N_{\e,i} =\int n_{\e,i} (z)dz=N_i^{*}(1+\e (F_i+\f{E_i}{A}+\f{3C}{A^2}))+O(\e^2),
\\
\mu_{\e,i} =\f{1}{N_{\e,i} }\int z  n_{\e,i} (z) dz=z^{*}+\e(\f{3B}{A^2}+\f{D_i}{A} )+O(\e^2),
\\
\sigma_{\e,i}^{  2}=\f{1}{N_{\e,i} }\int (z-\mu_{\e,i} )^2  n_{\e,i} (z) dz =\f{\e}{A} +O(\e^2),
\\
s_{\e,i} =\f{1}{\sigma_{\e,i}^{  3}N_{\e,i} }\int (z-\mu_{\e,i} )^3  n_{\e,i} (z) dz =\f{6B}{A^{\f32}}\sqrt{\e}+O(\e^{\f32}).
\end{cases}
$$
One can obtain similar approximations in the case of dimorphic ESS. To compute the above integrals, replacing the approximation \fer{ap-ue} in the integrals, a natural change of variable is to take $z-z^{*}=\sqrt{\e}y$. Therefore each term $z-z^{*}$ can be considered as of order $\sqrt{\e}$ in the integration. This is why, to obtain a first order approximation of the integrals in terms of $\e$, it is enough to have a fourth order approximation of $u(z)$, a second order approximation of $v_i(z)$ and a zero order approximation of $w_i(z)$,  in terms of $z$ around $z^{*}$.

% The Acknowledgements are an un-numbered section
\section*{Acknowledgements}
% Acknowledgements text here
S. Mirrahimi  has received funding from the European Research Council (ERC) under the European Union's Horizon 2020 research and innovation programme (grant agreement No 639638), and from the french ANR projects KIBORD ANR-13-BS01-0004 and MODEVOL ANR-13-JS01-0009.

%\bibliography{\string~/owncloud/Documents/aaaa-tex-files/Bibtex/Bibtex3/bibli}
%\bibliographystyle{plain}
% 

\end{document}